\documentclass[12pt]{article}
\usepackage{amsmath,amssymb}
\usepackage[mathscr]{euscript}
\usepackage{graphicx}
\usepackage[OT1,T1]{fontenc}
\usepackage[all]{xypic}
\usepackage{oubraces}
\usepackage{hyperref}
\usepackage{cancel}
\usepackage[normalem]{ulem}

\usepackage{xcolor}

\title{Dynamical mechanisms for Kaluza--Klein theories}
\author{Fr{\'e}d{\'e}ric \textsc{H{\'e}lein}\footnote{Institut de Math{\'e}matiques de Jussieu,
UMR CNRS 7586 Universit{\'e} Paris Cit{\'e},
UFR de Math{\'e}matiques,  B{\^a}timent Sophie Germain
75205 Paris Cedex 13, France, \textsf{helein@math.univ-paris-diderot.fr}}
}

\begin{document}
\maketitle

\emph{Abstract: We present variational formulations of gauge theories and Einstein--Yang--Mills equations in the
spirit of Kaluza--Klein theories. For gauge theories only a topological fibration is assumed. For
gravitation coupled with gauge fields no fibration is assumed: Fields are defined on a ’space-time’ $\mathcal{Y}$
of dimension $4 + r$ without any structure a priori, where $r$ is the dimension of the structure group.
If the latter is compact and simply connected, classical solutions allow to construct a manifold $\mathcal{X}$ of
dimension 4 to be the physical space-time, in such a way that $\mathcal{Y}$ acquires the structure of a principal
bundle over $\mathcal{X}$ and leads to solutions of the Einstein--Yang--Mills systems. The special case of the
Einstein-Maxwell system is also discussed : It suffices that at least one fiber closes in on a circle to
deduce that the five-dimensional space-time has a fiber bundle structure.}




\maketitle

\section{Introduction}\label{sec1}

Kaluza--Klein theories (see \cite{appelquist} for an account) go back to the work of T. Kaluza \cite{kaluza} in 1921 and O. Klein \cite{klein} in 1926, in an attempt to unify gravity and electromagnetism as a manifestation of a five-dimensional relativistic gravity theory on a manifold fibered over some four-dimensional space-time.
This theory was abandoned for various reasons before arousing new interest in the context of supergravity and string theory.
Some inconsistency was observed and fixed through the introduction an additional fields (\emph{radion} or \emph{dilaton}) independently by P. Jordan \cite{jordan} in 1947 and Y. Thiry \cite{thiry} in 1948. The introduction of this field, the physical meaning of which is not clear, can be avoided by not imposing the Einstein equation on the higher dimensional space-time and instead by looking for the critical points of the Einstein--Hilbert action when assuming a fiber bundle structure and equivariance constraints. By following this alternative option the theory was extended to Yang--Mills fields by R. Kerner \cite{kerner} in 1968, leading to the exact Einstein--Yang--Mills system (see also \cite{bleeker}).

Another problem is to explain why the extra dimensions cannot be observed. One of the most common hypothesis, due to Klein, is to suppose that the extra dimensions are tiny, so that, in the quantized version, higher modes are too energetic to be observable. Another one is simply to assume that the metric is invariant along extra dimensions. But this assumption needs to be justified.

This is the question we address here:
We propose theoretical models coupling gravity with electromagnetism or gauge fields such that, at the classical level, \emph{the invariance of fields along extra dimensions follows from dynamical equations}, so that there is no need to assume this condition \emph{a priori}.

Hence in the following models, fields which are not solutions of the dynamical equations depend in general on $4+r$ variables, where $r$ is the number of extra dimensions. But some of the dynamical equations, caused by \emph{auxiliary fields} which play
the role of Lagrange multipliers, force the apparition of a fibration of the higher dimensional space-time over some four-dimensional physical space-time.

However these auxiliary fields may spoil the result (in a way similar to \emph{dilatons}) by creating some artificial matter sources in the final equations. This is where an unexpected \emph{cancellation phenomenon} (\ref{cancel}) comes into play, leading to the vanishing of terms coming from auxiliary fields in the equations, when these are projected on the space-time.

\section{Results}\label{sec2}
In the following, we present several models of growing generality, in order to expound the various mechanisms underlying our theory. In Section \ref{deux} we present a model leading to Maxwell equations on a \emph{fixed} space-time but by replacing the usual principal $U(1)$-bundle by a topological circle bundle.
As a consequence the gauge symmetry group is larger than in the standard theory since, in addition to the usual gauge symmetries, the action functional is invariant by fibration preserving diffeomorphisms. This shows how the invariance of some fields along extra dimensions can be achieved and also gives a first glance at the \emph{cancellation mechanism}.

In Section \ref{trois} we extend this construction for Yang--Mills fields, still over a fixed space-time, for a structure group which is \emph{compact} and \emph{simply connected} (e.g., any product of $SU(k)$'s). It generalizes results in \cite{helein14} (where we needed to assume a principal bundle structure \emph{a priori}). 

In Section \ref{quatre} we then turn on the general Kaluza--Klein model for a \emph{compact} and \emph{simply connected} structure group. Fibers are obtained by integrating an exterior differential system through Frobenius's theorem. Most of these results were obtained in \cite{helein2020} (without a bare cosmological constant).

The previous model, however, does not work for $U(1)$ since this group is not simply connected.
This question is addressed in Section \ref{cinq} where we show that, if there exists \emph{at least} one integral curve of the exterior differential system which closes, then we obtain a four-dimensional space-time, a fibration and, at the very end, a solution of the Einstein--Maxwell system.

Throughout the paper all fields are assumed to be smooth.

\section{Model for Maxwell equations}\label{deux}
Consider a connected four-dimensional space-time manifold $\mathcal{X}^4$ with a fixed Riemannian metric $\mathbf{g}_{\mu\nu}$, a five-dimensional space-time $\mathcal{Y}^5$ and suppose that there exists a smooth projection map $\mathcal{Y}^5\xrightarrow{P}\mathcal{X}^4$
such that $\hbox{d}P$ has a maximal rank (equal to 4) everywhere and that, for any $x\in \mathcal{X}^4$ the inverse image $\textsf{f}_x:=P^{-1}(x)$ is a topological circle. We denote by $x = (x^\mu)_{0\leq \mu\leq 3}$ the coordinates on $\mathcal{X}^4$, by $y$ the coordinate on fibers $\textsf{f}_x$ and set $z = (z^I)_{0\leq I\leq 4} = (x,y)$. We describe the metric on $\mathcal{X}^4$ through a fixed co-vierbein $e^a = e^a_\mu(x) \emph{d}x^\mu$ for $0\leq a\leq 3$, so that $\mathbf{g}_{\mu\nu} = \eta_{ab}e^a_\mu e^b_\nu$, where $\eta_{ab}$ is the Minkowski metric.

Dynamical fields on $\mathcal{Y}^5$ are pairs $(\theta,\pi)$, where 
$\theta= \theta_\mu(x,y)\hbox{d}x^\mu + \theta_4(x,y)\hbox{d}y = \theta_I(z)\hbox{d}z^I$ and $\pi = \frac{1}{3!}\pi_{IJK}(z)\hbox{d}z^I\wedge \hbox{d}z^J\wedge \hbox{d}z^K$. We will, however, privilege another decomposition of $\pi$.
Indeed, we assume that $\theta_4\neq 0$ so that $(e_0,e_1,e_2,e_3,\theta)$ provides us with a co-5-bein on $\mathcal{Y}^5$.
We define
\[
 e^{(4)} = e^0\wedge e^1\wedge e^2\wedge e^3,
 \quad e_a^{(3)} = \frac{1}{3!}\epsilon_{abcd}e^{b}\wedge e^{c}\wedge e^{d}\quad
\hbox{ and }e_{ab}^{(2)} = \frac{1}{2!}\epsilon_{abcd}e^{c}\wedge e^{d}
\]
This allows us to decompose
\begin{equation}\label{eq:piMaxwell}
\begin{array}{ccc}
 \pi & = & \frac{1}{2}\pi^{ab}(x,y)e_{ab}^{(2)}\wedge \theta
-\pi^a(x,y)e_a^{(3)}
\end{array}
\end{equation}
We set $\|\pi\|^2 = \frac{1}{2}\pi^{ab}\pi_{ab}$, where $\pi_{ab} = \eta_{ac}\eta_{bd}\pi^{cd}$, and define the action functional 
\begin{equation}
\label{eq:actionMaxwell}
 \mathcal{A}[\theta,\pi]
= \int_{\mathcal{Y}^5}
 \frac{1}{2}\|\pi\|^2e^{(4)}\wedge \theta + \pi\wedge \hbox{d}\theta
\end{equation}
or $\mathcal{A}[\theta,\pi] = \int_{\mathcal{Y}^5}\left(
 \frac{1}{4}\pi^{ab}\pi_{ab} + \frac{1}{2}\Theta_{ab}\pi^{ab} + \Theta_a\pi^a \right)e^{(4)}\wedge \theta$,
where we set 
$\hbox{d}\theta
= \frac{1}{2}\Theta_{ab}e^a\wedge e^b + \Theta_a e^a\wedge \theta$.\\

\subsection{Euler--Lagrange equations}
We denote by $(e_I)_{0\leq I\leq 4}$ the 5-bein dual to $(e^a,\theta)$. The Euler--Lagrange equations are
\begin{equation}\label{ELMaxwell}
 \begin{array}{rcccc}
  \hbox{d}\theta(e_a,\partial_y)/\theta_4 = 
  \Theta_a
  & = & 0 & \hbox{(a)} \\
\Theta_{ab} + \pi_{ab}
  & = & 0 & \hbox{(b)} \\
\hbox{d}\pi & = & \frac{1}{2}\|\pi\|^2e^{(4)} & \hbox{(c)}
 \end{array}
\end{equation}
where Equations (a), (b) and (c) correspond to variations with respect to $\pi_a$, $\pi_{ab}$ and $\theta$, respectively. 
Equation (\ref{ELMaxwell}a) has the following consequence. Consider two points $x_1,x_2\in \mathcal{X}^4$ and a path $\gamma$ joining $x_1$ to $x_2$ in $\mathcal{X}^4$. Its inverse image by $P$ is a surface $S\subset \mathcal{Y}^5$ with boundary $\partial S = \textsf{f}_{x_2} - \textsf{f}_{x_1}$. Since $\partial_y$ is tangent to $S$, we have by (\ref{ELMaxwell}a)
\begin{equation}\label{fluxMaxwell}
 \int_{\textsf{f}_{x_2}}\theta -
 \int_{\textsf{f}_{x_1}}\theta
 = \int_S\hbox{d}\theta = 0
\end{equation}
Hence, $q:= \int_{\textsf{f}_x}\theta$ is independent of $x$ in $\mathcal{X}^4$. We can thus define the new coordinate
\[
 s = f(x,y) = \int_0^y \theta_4(x,y')\hbox{d}y'
 \hbox{ mod }[q]
\]
and replace the coordinates $(x,y)$ by $(x,s)$, where $s\in \mathbb{R}/q\mathbb{Z}$. The form $\theta$ then reads $\theta = \mathbf{A}_a(x,s)e^a + \hbox{d}s$, where $\textbf{A}_a(x,f(x,y)) = (\theta- \hbox{d}f)_{(x,y)}(e_a)$.
Equation (\ref{ELMaxwell}a) further implies $\partial_s\mathbf{A}_a = 0$, so that 
\begin{equation}\label{Maxwellnewtheta}
 \theta = \mathbf{A}_a(x)e^a + \hbox{d}s
\end{equation}
Hence $\hbox{d}\theta = \mathbf{F} = \frac{1}{2}\mathbf{F}_{ab}e^a\wedge e^b$, where $\mathbf{F}_{ab} = \partial_a\mathbf{A}_b - \partial_a\mathbf{A}_b$. We define $p^{ab}(x,s)$ and $p^a(x,s)$ such that
\[
 \begin{array}{ccl}
  p^a(x,f(x,y)) & = & \pi^a(x,y)-\pi^{ab}(x,y)\textbf{A}_b(x,f(x,y)) \\
 p^{ab}(x,f(x,y)) & = & \pi^{ab}(x,y)
 \end{array}
\]
Then $\pi = \frac{1}{2}p^{ab}e_{ab}^{(2)}\wedge \hbox{d}s - p^ae_a^{(3)}$, where $p^{ab} = -\mathbf{F}^{ab}$ by (\ref{ELMaxwell}b).

Let $\gamma{^a}_b$ be the spin connection on $\mathcal{X}^4$ for the vierbein $(e_a)$.
It satisfies the torsion free condition
\begin{equation}\label{torsionfree}
\hbox{d}^\gamma e^a:= \hbox{d} e^a + \gamma{^a}_b\wedge e^b = 0.
\end{equation}
By setting 
$\hbox{d}^\gamma p^{ab}:= \hbox{d}p^{ab} + \gamma{^a}_cp^{cb} + \gamma{^a}_cp^{ac}$ and $\hbox{d}^\gamma p^a:= \hbox{d}p^a + \gamma{^a}_cp^c$ and by using (\ref{torsionfree}) we find that 
$\hbox{d}\pi = \frac{1}{2}\hbox{d}^\gamma p^{ab}\wedge e_{ab}^{(2)}\wedge \hbox{d}s - \hbox{d}^\gamma p^a\wedge e_a^{(3)}$.
By decomposing
$\hbox{d}^\gamma p^{ab} = \partial_c^\gamma p^{ab}e^c + \partial_s p^{ab}\hbox{d}s$
and  
$\hbox{d}^\gamma p^a = \partial_c^\gamma p^ae^c + \partial_s p^a\hbox{d}s$ we 
hence get 
$\hbox{d}\pi = (\partial_b^\gamma p^{ab}+\partial_sp^a) e_{a}^{(3)}\wedge \hbox{d}s - \partial_a^\gamma p^a e^{(4)}$. Thus, (\ref{ELMaxwell}c) is equivalent to
\begin{equation}\label{ELMaxplus}
 \begin{array}{cccc}
  \partial_b^\gamma \textbf{F}^{ab}(x) & = & \partial_sp^a(x,s) & \hbox{(a)} \\
    \partial_a^\gamma p^a(x,s) & = & -\frac{1}{2}\|\textbf{F}\|^2(x) & \hbox{(b)}
 \end{array}
\end{equation}
The key point to conclude that $\mathbf{A}$ is a solution of Maxwell equations in vacuum is to observe that the l.h.s of (\ref{ELMaxplus}a) is independent of $s$, thus
\begin{equation}
 \partial_b^\gamma \textbf{F}^{ab}
 = \frac{\int_0^q\partial_b^\gamma \textbf{F}^{ab}\hbox{d}s}{\int_0^q\hbox{d}s}
 = \frac{\int_0^q\partial_sp^a \hbox{d}s}{q}
 = 0
\end{equation}
because $p^a$ is $q$-periodic in $s$. Equation (\ref{ELMaxplus}b) involves fields which are not observable \emph{a priori}.\\

\subsection{Gauge symmetries}
The action (\ref{eq:actionMaxwell}) is invariant under several types of gauge symmetries:

(i) by diffeomorphisms of $\mathcal{Y}^5$ which preserves the orientation and the map $\mathcal{Y}^5\xrightarrow{P}\mathcal{X}^4$, i.e. of the form $T(x,y) = (x,f(x,y))$, through pullback $(\theta,\pi)\longmapsto (T^*\theta,T^*\pi)$;

(ii) by gauge transformations $\theta\longmapsto \theta + \hbox{d}V$, where $V$ is a function of $x\in \mathcal{X}^4$;

(iii) by transformations $\pi\longmapsto \pi + \psi$, where $\psi$ is any closed 3-form on $\mathcal{X}^4$ which decays at infinity.

\section{Non-Abelian gauge fields}\label{trois}
The previous theory can be extended to nonlinear gauge theories. We let $\mathfrak{G}$ be a compact simply connected Lie group of dimension $r$ and we replace $\mathcal{Y}^5
$ by a manifold $\mathcal{Y}^{N+1}$ of dimension $N+1 =4+r$. We still assume the existence of a smooth projection map $\mathcal{Y}^{N+1}\xrightarrow{P}\mathcal{X}^4$ such that $\hbox{d}P$ has maximal rank and, for any $x\in \mathcal{X}^4$, $\textsf{f}_x:= P^{-1}(x)$ is a \emph{connected} submanifold of $\mathcal{Y}^{N+1}$ of dimension $r$. Local coordinates on $\mathcal{Y}^{N+1}$ are $z = (z^I)_{0\leq I\leq N} = (x,y) = (x^\mu,y^i)_{0\leq \mu\leq 3< i\leq N}$, where $(x^\mu)_{0\leq \mu\leq 3}$ are coordinates on $\mathcal{X}^4$ as previously.

Let $\mathfrak{g}$ be the Lie algebra of $\mathfrak{G}$, $(\mathbf{t}_i)_{3< i\leq N}$ a basis of $\mathfrak{g}$ and $c^i_{jk}$ the structure coefficients such that $[\mathbf{t}_i,\mathbf{t}_j] = c^k_{ij}\mathbf{t}_k$.  We note $\mathfrak{g}^*$ the dual space of $\mathfrak{g}$ and $(\mathbf{t}^i)_{3< i\leq N}$ its dual basis. A consequence of the compactness of $\mathfrak{G}$ is that $\mathfrak{g}$ is \emph{unimodular}, which reads
\begin{equation}\label{unimodular}
c^i_{ij} = c^i_{ji} = 0
\end{equation}
We let $\textsf{k}$ be a scalar product on $\mathfrak{g}$ which invariant by the adjoint action of $\mathfrak{G}$
and we set $\textsf{k}_{ij} = \textsf{k}(\mathbf{t}_i,\mathbf{t}_j)$.
Fields are pairs $(\theta,\pi)$ where
\[
\begin{array}{ccl}
 \theta & = & \theta^i \mathbf{t}_i
 = (\theta{^i}_\mu(z)\hbox{d}x^\mu + \theta{^i}_j(z)\hbox{d}y^j)\mathbf{t}_i \\
 \pi & = & \pi_i\mathbf{t}^i
 = \frac{1}{(N-1)!}\pi_{iI_1\cdots I_{N-1}}(z)\hbox{d}z^{I_1}\wedge \cdots\wedge \hbox{d}z^{I_{N-1}}\mathbf{t}^i
\end{array}
\]
As for Maxwell case we assume that the rank of $(e^a,\theta^i)_{0\leq a\leq 3< i\leq N}$ is maximal, equal to $N+1$. Hence by introducing the notations
\begin{equation}\label{conventions}
 \bar{\theta}^{(r-\alpha)}_{i_1\cdots i_\alpha}
 = \frac{1}{(r-\alpha)!}\epsilon_{i_1\cdots i_\alpha j_{\alpha +1} \cdots j_r}
 \theta^{j_{\alpha+1}}\wedge \cdots \wedge \theta^{j_r}
\end{equation}
for $0\leq \alpha \leq r$ and where all indices run from 4 to $N$, 
we instead decompose $\pi$ as
\begin{equation}
 \pi_i = \frac{1}{2} \pi{_i}^{ab}e_{ab}^{(2)}\wedge \bar{\theta}^{(r)}
 - \pi{_i}^{ak}e_a^{(3)}\wedge \bar{\theta}^{(r-1)}_k + \frac{1}{2}\pi{_i}^{jk}e^{(4)}\wedge \bar{\theta}^{(r-2)}_{jk}
\end{equation}
where coefficients $\pi{_i}^{JK}$ are functions of $z = (x,y)$.

The action functional is
\begin{equation}\label{ActionYM}
 \mathcal{A}[\theta,\pi] = \int_{\mathcal{Y}^N}\frac{1}{2}\|\pi\|^2e^{(4)}\wedge \bar{\theta}^{(r)}
 + \pi_i\wedge \Theta^i
\end{equation}
where $\|\pi\|^2 = \frac{1}{2}\pi{_i}^{ab}\pi{^i}_{ab}$
with $\pi{^i}_{ab} = \textsf{k}^{ij}\pi{_j}^{cd}\eta_{ac}\eta_{bd}$
and $\Theta^i:= \hbox{d}\theta^i+\frac{1}{2}[\theta\wedge\theta]^i = 
\hbox{d}\theta^i+\frac{1}{2}c^i_{jk}\theta^j\wedge\theta^k$.\\

\subsection{Euler--Lagrange equations}
Using the decomposition 
$\Theta^i = \frac{1}{2}\Theta{^i}_{ab}e^a\wedge e^b + \Theta{^i}_{ak}e^a\wedge \theta^k + \frac{1}{2}\Theta{^i}_{jk}\theta^j\wedge \theta^k$ and denoting $\hbox{d}^\theta \pi_i
= \hbox{d}\pi_i - c^j_{ki}\theta^k\wedge \pi_j$, the dynamical equations read
\begin{equation}\label{ELYM}
 \begin{array}{cccc}
\pi{^i}_{ab} + \Theta{^i}_{ab} & = & 0 & \hbox{(a)} \\
\Theta{^i}_{ak} & = & 0 & \hbox{(b)} \\
\Theta{^i}_{jk} & = & 0 & \hbox{(c)} \\
\hbox{d}^\theta \pi_i & = & \frac{1}{2}\|\pi\|^2e^{(4)}\wedge \bar{\theta}^{(r-1)}_i & \hbox{(d)}
 \end{array}
\end{equation}
where Equations (a), (b), (c) and (d) correspond to variations with respect to $\pi{_i}^{ab}$, $\pi{_i}^{ak}$, $\pi{_i}^{jk}$ and $\theta^i$, respectively.
We first use (\ref{ELYM}c), which implies that, for any $3<i\leq N$ and for any $x\in \mathcal{X}^4$, the restriction of $\Theta^i = \hbox{d}\theta^i + \frac{1}{2}c^i_{jk}\theta^j\wedge \theta^k$ on $\textsf{f} = \textsf{f}_x$ vanishes. Hence, by Frobenius' theorem (see \cite{bcggg}), given a point $z\in \textsf{f}$, we can construct a unique map $g$ from a neighborhood of $z$ in $\textsf{f}$ to $\mathfrak{G}$ 
such that the restriction to $\textsf{f}$ of
$\theta - g^{-1}\hbox{d}g$ vanishes
and such that $g(z) = 1$. Since the rank of $(\theta^i )_{3 < i \leq N}$ is $r$, this map is a local diffeomorphism.
Actually, the inverse map can be extended globally as a map $T_\textsf{f}:\mathfrak{G}\rightarrow \textsf{f}$ which associates with any $g\in \mathfrak{G}$ the end value $v(1)$ of the path $v:[0,1]\longrightarrow \textsf{f}$ which is a solution of $v(0) = z$, $e^a(\frac{dv}{dt}) = 0$ for $0\leq a\leq3$ and $\theta(\frac{dv}{dt}) = u^{-1}\frac{du}{dt}$ where $u:[0,1]\longrightarrow \mathfrak{G}$ is a path such that $u(0) = 1_\mathfrak{G}$ and $u(1) = g$. Indeed the definition of $T_\textsf{f}(g)$ does not depend on the choice of the path $u$ since $\mathfrak{G}$ is simply connected.
$T_\textsf{f}$ is then a covering map of $\textsf{f}$. Hence since $\mathfrak{G}$ is compact, $\textsf{f}$ is compact and is diffeomorphic to a quotient $\mathfrak{G}_\textsf{f}$ of $\mathfrak{G}$ by a finite subgroup.

But all fibers are diffeomorphic to the same group $\mathfrak{G}_0$. Indeed, for any fixed $\xi = \xi^a\mathbf{t}_a$ in $\mathfrak{g}$, consider the vector field $\textsf{X} = \textsf{X}(\xi)$ on $\mathcal{Y}^{N+1}$ defined by $e^a(\textsf{X})=\xi^a$ and $\theta^i(\textsf{X})= 0$. Then by (\ref{torsionfree}) which, by setting $\gamma{^a}_b = \gamma{^a}_{bc}e^c$, reads $\hbox{d}e^a = \gamma{^a}_{bc}e^b\wedge e^c$, the Lie derivative of $e^a$ by $X$ satisfies $L_\textsf{X}e^a =   0\hbox{ mod }[e^b]$, for any $0\leq a,b\leq 3$. This implies that the image of a fiber $\textsf{f}$ by the flow map of $\textsf{X}$ is another fiber. We can thus construct a diffeomorphism between $\mathfrak{G}_\textsf{f}\times B^4(0,\varepsilon)\simeq \textsf{f}\times B^4(0,\varepsilon)$ (where $B^4(0,\varepsilon)$ is the ball of a sufficiently small radius $\varepsilon$ in $\mathbb{R}^4$) and a neighborhood of a leaf $\textsf{f}$ in $\mathcal{Y}^{N+1}$ by mapping $(z,\xi)$ to $e^{\textsf{X}(\xi)}(z)$. Its inverse map provides us with a local trivialization of $\mathcal{Y}^{N+1}$. Since this construction can be done everywhere, $\mathcal{Y}^{N+1}\rightarrow\mathcal{X}^4$ is endowed with a structure of principal bundle with a structure group $\mathfrak{G}_0$.

By choosing a local trivialization of this bundle, we can improve the use of Equation (\ref{ELYM}c) to show the existence of a map $g$ from an open subset of $\mathcal{Y}^{N+1}$ to $\mathfrak{G}_0$ such that the restriction of $\theta -g^{-1}\hbox{d}g$ to $\textsf{f}$ vanishes for any fiber $\textsf{f}$.
This means that, if we set $\mathbf{A} = g\theta g^{-1} - \hbox{d}g\ g^{-1}$ and decompose $\mathbf{A} = \mathbf{A}_ae^a + \mathbf{A}_i\theta^i$, then $\mathbf{A}_i=0$. Hence
$\theta_z = g^{-1}(z)\mathbf{A}_a(z) g(z)e^a + g^{-1}(z)\hbox{d}g_z$.

Equation (\ref{ELYM}b) then reads $\partial_k\mathbf{A}_a = 0$, so that actually $\mathbf{A}_a(z) = \mathbf{A}_a(x)$ and 
\begin{equation}
 \theta = g^{-1}\mathbf{A}_a(x)ge^a + g^{-1}\hbox{d}g.
\end{equation}
Thus, $\mathcal{Y}^{N+1}\rightarrow \mathcal{X}^4$ is endowed with the connection form $\mathbf{A}$. We define $\mathbf{F} = \hbox{d}\mathbf{A} + \frac{1}{2}[\mathbf{A}\wedge \mathbf{A}] = \frac{1}{2}\mathbf{F}_{ab}(x)e^a\wedge e^b$, 
so that we have $\Theta:= \Theta^i\mathbf{t}_i = g^{-1}\mathbf{F}g$.

Let $\left(S^i_j\right)$  be the matrix of the adjoint action of $g$ on $\mathfrak{g}$ in the basis $(\mathbf{t}_i)_{3<i\leq N}$, i.e., such that $g\mathbf{t}_jg^{-1} = S^i_j\mathbf{t}_i$. We define $e^i = S^i_je^j$ and, by using the same conventions as in (\ref{conventions}), $\bar{e}^{(r)}$, 
$\bar{e}^{(r-1)}_i$ and $\bar{e}^{(r-2)}_{ij}$.
Then, by (\ref{unimodular}), $\bar{e}^{(r)} = \bar{\theta}^{(r)}$, $\bar{e}^{(r-1)}_i=(S^{-1})_i^j\bar{\theta}^{(r-1)}_j$ and $\bar{e}^{(r-2)}_{ij}=(S^{-1})_i^k(S^{-1})_j^\ell\bar{\theta}^{(r-2)}_{k\ell}$. We also set
\[
 \begin{array}{llllll}
 p_i & = & (S^{-1})^j_i\pi_j &  p{_i}^{ab} & = & (S^{-1})^j_i\pi{_j}^{ab}\\
 p{_i}^{ak} & = & (S^{-1})^j_iS^k_\ell \pi{_j}^{a\ell} \quad & 
 p{_i}^{jk} & = & (S^{-1})^{\ell}_iS^j_{m}S^k_{n} \pi{_\ell}^{mn}
  \quad 
 \end{array}
\]
Then,
\begin{equation}\label{piYM}
 p_i\! =\! \frac{1}{2}p{_i}^{ab}e^{(2)}_{ab}\wedge \bar{e}^{(r)} - p{_i}^{ak}e^{(3)}_{a}\wedge \bar{e}^{(r-1)}_k
 + \frac{1}{2}p{_i}^{jk}e^{(4)}\wedge \bar{e}^{(r-2)}_{jk}
\end{equation}
Then, (\ref{ELYM}a) reads $p{^i}_{ab} + \mathbf{F}{^i}_{ab} = 0$, where $p{^i}_{ab} = \textsf{k}^{ij}\eta_{ac}\eta_{bd}p{_j}^{cd}$,  or, by setting $\mathbf{F}{_i}^{ab} = \textsf{k}_{ij}\eta^{ac}\eta^{bd}\mathbf{F}{^j}_{cd}$,
\begin{equation}\label{YMpplusF}
 p{_i}^{ab} + \mathbf{F}{_i}^{ab} = 0 
\end{equation}
Lastly, by defining $\hbox{d}^\mathbf{A}p_i = \hbox{d}p_i -c^j_{ki}\mathbf{A}^k\wedge p_j$ and $\|p\|^2 = \frac{1}{2}p{^i}_{ab}p{_i}^{ab}$, Equation (\ref{ELYM}d) translates as
\begin{equation}\label{preYM}
 \hbox{d}^\mathbf{A}p_i = \frac{1}{2}\|p\|^2e^{(4)}\wedge \bar{e}_i^{(r-1)}
\end{equation}
The computation of $\hbox{d}^\mathbf{A}p_i$ requires some further notations. 
We denote by $\partial_a$ and $\partial_i$ the operators such that, for any function $f$, $\hbox{d}f = (\partial_af)e^a+ (\partial_if)e^i$ and by $\mathbf{A}{^i}_a$ the coefficients such that $\mathbf{A}^i = \mathbf{A}{^i}_ae^a$. As in Section \ref{deux} we let $\gamma{^a}_b$ be the spin connection coefficients and we set $\gamma{^a}_c = \gamma{^a}_{cb}e^b$. Lastly we define
%
\[
\begin{array}{lll}
  \partial^{\gamma,\mathbf{A}}_b p{_i}^{ab} & = &  \partial_bp{_i}^{ab} - \mathbf{A}{^k}_bc^\ell_{ki}p{_\ell}^{ab}  + \gamma{^a}_{cb}p{_i}^{cb} + \gamma{^b}_{cb}p{_i}^{ac} \\
\partial^{\gamma,\mathbf{A}}_bp{_i}^{jb} & = & \partial_bp{_i}^{jb} - \mathbf{A}{^k}_bc^\ell_{ki}p{_\ell}^{jb} + \mathbf{A}{^k}_bc^j_{k\ell} p{_i}^{\ell b} + \gamma{^b}_{cb}p{_i}^{jc} 
\end{array}
\]
Then, we obtain by using (\ref{torsionfree})  and (\ref{unimodular}) that
\begin{equation}\label{15abis}
\begin{array}{ccc}
  \hbox{d}^\mathbf{A}p_i & = & (\partial^{\gamma,\mathbf{A}}_bp{_i}^{ab} + \partial_kp{_i}^{ak})\ e^{(3)}_a\wedge \bar{e}^{(r)}  \\
  & & +\ (\partial^{\gamma,\mathbf{A}}_bp{_i}^{jb}+ \partial_kp{_i}^{jk} + \frac{1}{2}\mathbf{F}{^j}_{ab}p{_i}^{ab}
  + \frac{1}{2}c^j_{k\ell}
  p{_i}^{k\ell})\ e^{(4)}\wedge \bar{e}^{(r-1)}_j
\end{array}
\end{equation}
By taking into account (\ref{YMpplusF}), we deduce that (\ref{preYM}) is equivalent to
\begin{equation}\label{YMfinalsystem}
\begin{array}{rll}
\partial^{\gamma,\mathbf{A}}_b\mathbf{F}{_i}^{ab} & = & \partial_kp{_i}^{ak} \\
\partial^{\gamma,\mathbf{A}}_bp{_i}^{jb} + \partial_kp{_i}^{jk} + \frac{1}{2}c^j_{k\ell}p{_i}^{k\ell} & = &  \frac{1}{2}\|\mathbf{F}\|^2\delta{_i}^j
+ \frac{1}{2}\mathbf{F}{_i}^{ab}\mathbf{F}{^j}_{ab} 
\end{array}
\end{equation}
As in Maxwell theory, a key point is to observe that the l.h.s. in the first equation in (\ref{YMfinalsystem}) does not depend of $y$ (or $g$). Hence, since the fibers are  compact,
\begin{equation}\label{cancel}
\partial^{\gamma,\mathbf{A}}_b\mathbf{F}{_i}^{ab} = \frac{\int_{\textsf{f}_x}\partial^{\gamma,\mathbf{A}}_b\mathbf{F}{_i}^{ab}\ \bar{e}^{(r)}}{\int_{\textsf{f}_x}\bar{e}^{(r)}} = \frac{\int_{\textsf{f}_x}\partial_kp{_i}^{ak}\ \bar{e}^{(r)}}{\int_{\textsf{f}_x}\bar{e}^{(r)}}  = \frac{\int_{\textsf{f}_x}\hbox{d}(p{_i}^{ak}\bar{e}^{(r-1)}_k)}{\int_{\textsf{f}_x}\bar{e}^{(r)}} = 0
\end{equation}
Hence, $\mathbf{A}$ is a solution of Yang--Mills equations. The second equation in (\ref{YMfinalsystem}) involves fields which are not observable.\\

\subsection{Gauge symmetries}
The action (\ref{ActionYM}) is invariant under the following gauge symmetries:

(i) by diffeomorphisms of $\mathcal{Y}^{N+1}$ which preserves the orientation and the map $\mathcal{Y}^{N+1}\xrightarrow{P}\mathcal{X}^4$, i.e., of the form $T(x,y) = (x,f(x,y))$, through pullback $(\theta,\pi)\longmapsto (T^*\theta,T^*\pi)$;

(ii) by gauge transformations 
\begin{equation}
\begin{array}{ccl}
\theta & \longmapsto & \hbox{Ad}_g\theta - \hbox{d}g\ g^{-1} \\
\pi & \longmapsto & \hbox{Ad}_g^*\pi
\end{array}
\end{equation}
for any map $g:\mathcal{X}^4\longrightarrow \mathfrak{G}$ and
where $\hbox{Ad}_g\theta = g\theta g^{-1}$ and, if $\hbox{Ad}_g\mathbf{t}_j = S^i_j\mathbf{t}_i$, $\hbox{Ad}_g^*\pi = (S^{-1})^j_i\pi_j\mathbf{t}^i$;

(iii) by transformations $\pi\longmapsto \pi + \chi$, where $\chi$ has the form
\[
\chi = \chi{_i}^{bk}e_b^{(3)}\wedge \bar{\theta}^{(r-1)}_k 
 + \frac{1}{2} \chi{_i}^{jk}e^{(4)}\wedge \bar{\theta}^{(r-2)}_{jk}
\]
decays at infinity and satisfies $\theta^i\wedge \hbox{d}^\theta \chi_i = 0$.

\section{Einstein--Yang--Mills model}\label{quatre}
The previous theory can be modified in order to couple gravitation with gauge fields. We consider the metric $\textsf{h}$ on $\mathbb{R}^4\times \mathfrak{g}$ such that $\mathbb{R}^4\perp \mathfrak{g}$, $\textsf{h}$ coincides with the Minkowski metric $\eta_{ab}$ on $\mathbb{R}^4$ and with $\textsf{k}$ on $\mathfrak{g}$.
We still assume here that $\mathfrak{G}$ is compact and simply connected, and we still work on a higher-dimensional `space-time' $\mathcal{Y}^{N+1}$. However:

(i) We replace the fixed co-vierbein $(e^a)_{0\leq a\leq 3}$ by a dynamical one $(\theta^a)_{0\leq a\leq 3}$, and we assume that the dynamical fields $(\theta^I)_{0\leq I\leq N} = (\theta^a,\theta^i)_{0\leq a\leq 3<i\leq N}$ forms a co-$(N+1)$-bein on $\mathcal{Y}^{N+1}$.

(ii) No more fibration of $\mathcal{Y}^{N+1}$ nor the existence of $\mathcal{X}^4$ is assumed \emph{a priori}. Even better, we do not make assumption on its topology, beside the fact that $\mathcal{Y}^{N+1}$ is oriented and connected. Instead, we introduce extra auxiliary fields $(\pi_a)_{0\leq a\leq 3}$ which are $(N-1)$-forms and which play the role of Lagrange multipliers for creating a fibration. 

(iii) We also introduce the field
 $\varphi{^I}_J = \varphi{^{I}}_J(z)\hbox{d}z^K$ with coefficients in the Lie algebra $so(\mathbb{R}^4\oplus \mathfrak{g},\textsf{h})$ of the Lorentz group, i.e. such that $\varphi^{IJ} + \varphi^{JI} = 0$, where
 $\varphi{^I}_J = \varphi^{IK}\textsf{h}_{KJ}$.
 
(iv) In the action we replace the term $\frac{1}{2}\|\pi\|^2e^{(4)}\wedge
  \bar{\theta}^{(r)}$ by the $(N+1)$-dimensional \emph{Einstein--Palatini} density
$\frac{1}{2}\hat{\theta}^{(N-1)}_{IJ}\wedge \Phi^{IJ}$, where $\Phi^{IJ}:= \hbox{d}\varphi^{IJ} + \varphi{^I}_K\wedge \varphi^{KJ}$ and $\hat{\theta}^{(N+1)} = \theta^0\wedge \cdots \wedge \theta^N$ and, for $0\leq \alpha$,
\[
 \hat{\theta}^{(N-\alpha)}_{I_0\cdots I_\alpha} = \frac{1}{(N-\alpha)!}\epsilon_{I_0\cdots I_\alpha J_{\alpha+1}\cdots J_N}
 \theta^{J_{\alpha+1}}\wedge \cdots \wedge \theta^{J_N}
\]
(v) Lastly, we impose the constraint
  \begin{equation}\label{constraint}
   \theta^a\wedge \theta^b\wedge \pi_c
   = \theta^a\wedge \theta^b\wedge \pi_i = 0
  \end{equation}

\emph{To summarize}: dynamical fields are the 1-forms
$(\theta^I)_{0\leq I\leq N} = (\theta^a,\theta^i)_{0\leq a\leq 3<i\leq N}$, of rank $N+1$ everywhere, $(\varphi^{IJ})_{0\leq I,J\leq N}$ with $\varphi^{IJ} + \varphi^{JI} = 0$ and the $(N-1)$-forms $(\pi_I)_{0\leq I\leq N} = (\pi_a,\pi_i)_{0\leq a\leq 3<i\leq N}$. Constraint (\ref{constraint}) is assumed, i.e., $\theta^a\wedge \theta^b\wedge \pi_I=0$.
The basic action is
\begin{equation}\label{actionfinale}
 \mathcal{A}_0[\theta,\varphi,\pi]
 = \int_{\mathcal{Y}^{N+1}}
\frac{1}{2}\hat{\theta}^{(N-1)}_{IJ}\wedge \Phi^{IJ} + \pi_I\wedge \Theta^I
\end{equation}
where
\[
 \begin{array}{cclc}
  \Theta^a & := & \hbox{d}\theta^a & \hbox{for }0\leq a\leq 3 \\
  \Theta^i & := & \hbox{d}\theta^i 
 +\frac{1}{2}c^i_{jk}\theta^j\wedge \theta^k & \hbox{for }3<i\leq N
 \end{array}
\]
We shall also incorporate a bare cosmological constant $\Lambda_0$ and consider the action
\begin{equation}\label{22bis}
  \mathcal{A}_{\Lambda_0}[\theta,\varphi,\pi]
 = \mathcal{A}_0[\theta,\varphi,\pi]
 - \int_{\mathcal{Y}^{N+1}}
 \Lambda_0 \hat{\theta}^{(N+1)}
\end{equation}
Constraint (\ref{constraint}) can actually be rephrased through the decomposition
\begin{equation}\label{25}
 \pi_I = - \pi{_I}^{ak}\theta^{(3)}_a\wedge \bar{\theta}^{(N)}_k + \frac{1}{2}\pi{_I}^{jk}\theta^{(4)}\wedge \bar{\theta}^{(N-1)}_{jk}
\end{equation}
(i.e., $\pi{_I}^{ab}$ vanishes).
We also use the decomposition
\begin{equation}
 \Theta^I = \frac{1}{2}\Theta{^I}_{JK}\theta^J\wedge \theta^K
\end{equation}
or $\Theta^I = \frac{1}{2}\Theta{^I}_{ab}\theta^a\wedge \theta^b
 + \Theta{^I}_{ak}\theta^a\wedge \theta^k
 + \frac{1}{2}\Theta{^I}_{jk}\theta^j\wedge \theta^k$.\\

\subsection{Euler--Lagrange equations}
We set 
$\hbox{d}^\varphi \hat{\theta}^{(N-1)}_{IJ}:= \hbox{d}\hat{\theta}^{(N-1)}_{IJ} - \varphi{^K}_I\wedge \hat{\theta}^{(N-1)}_{KJ}
- \varphi{^K}_J\wedge \hat{\theta}^{(N-1)}_{IK}$, $\hbox{d}^\theta \pi_a = \hbox{d} \pi_a$ if $0\leq a\leq 3$, $\hbox{d}^\theta \pi_i = \hbox{d}\pi_i - c^j_{ki}\theta^k\wedge \pi_j$ if $3< i\leq N$.
Then, the Euler--Lagrange equations of $\mathcal{A}_{\Lambda_0}[\theta,\varphi,\pi]$ read:
 \begin{equation}\label{ELEYM}
 \begin{array}{rccclcccc}
\Theta{^I}_{ak} & = &  0 & \hbox{(a)}
\\
\Theta{^I}_{jk} & = & 0 & \hbox{(b)} 
\\
\hbox{d}^\varphi \hat{\theta}^{(N-1)}_{IJ} & = & 0 & \hbox{(c)}
\\
\hbox{d}^\theta \pi_I + \frac{1}{2}\hat{\theta}^{(N-2)}_{IJK}\wedge\Phi^{JK} - \Theta{^K}_{IJ}\ \pi{_K}^{\ell J}\ \hat{\theta}^{(N)}_\ell
& = & \Lambda_0\hat{\theta}^{(N)}_I & \hbox{(d)}
 \end{array}
\end{equation}
where Equations (a), (b), (c) and (d) correspond to variations with respect to $\pi{_I}^{ak}$, $\pi{_I}^{jk}$, $\varphi^{IJ}$ and $\theta^I$, respectively.
Consider a solution of this system.
We endow $\mathcal{Y}^{N+1}$ with the metric $\mathbf{h}_{\mu\nu}\hbox{d}z^\mu\otimes\hbox{d}z^\nu:= \textsf{h}_{IJ}\theta{^I}\otimes\theta{^J}$, and we assume that $(\mathcal{Y}^{N+1},\mathbf{h})$ is \emph{complete}, meaning that any geodesic curve is defined for all `time'.

Equation (\ref{ELEYM}c) is equivalent to $\hbox{d}^\varphi\theta^I:=\hbox{d}\theta^I + \varphi{^I}_J\wedge \theta^J = 0$ thanks to the identity
$\hbox{d}^\varphi \hat{\theta}^{(N-1)}_{IJ} = \hbox{d}^\varphi\theta^K\wedge \hat{\theta}^{(N-2)}_{IJK}$. It means that the spin connection defined by $\varphi^{IJ}$ is torsionfree and, hence, coincides with the Levi-Civita connection on $(\mathcal{Y}^{N+1},\mathbf{h})$.

Equation (\ref{ELEYM}b) implies first $\hbox{d}\theta^a = 0\hbox{ mod }[\theta^b]$, for any $0\leq a,b\leq 3$. This allows us to apply Frobenius' theorem (see \cite{bcggg}) to prove that $\mathcal{Y}^{N+1}$ is foliated by $r$-dimensional submanifolds (leaves) $\textsf{f}$ such that the restriction of $\theta^a$ to $\textsf{f}$ vanishes, for any $0\leq a\leq 3$. We then define $\mathcal{X}$ to be the set of leaves of this foliation (not yet a manifold !).

By using the same reasoning as in Section \ref{trois} (replacing (\ref{ELYM}c) by (\ref{ELEYM}b), which implies that the restriction to $\textsf{f}$ of $\Theta^i = \hbox{d}\theta^i + \frac{1}{2}c^i_{jk}\theta^j\wedge \theta^k$ vanishes) we prove that each leaf $\textsf{f}$ is compact and diffeomorphic to a quotient $\mathfrak{G}_\textsf{f}$ of $\mathfrak{G}$ by a finite subgroup.

We then prove that the foliation forms actually a fibration. For any fixed $\xi = \xi^a\mathbf{t}_a$, consider the vector field $\textsf{X} = \textsf{X}(\xi)$ on $\mathcal{Y}^{N+1}$ defined by $\theta^a(\textsf{X})=\xi^a$ and $\theta^i(\textsf{X})= 0$. By (\ref{ELEYM}a,b), we have $\hbox{d}\theta^a = \frac{1}{2}\Theta{^a}_{bc}\theta^b\wedge \theta^c$ and thus $L_\textsf{X}\theta^a = \Theta{^a}_{bc}\xi^b\theta^c$. Hence,
$L_\textsf{X}\theta^a = 0 \hbox{ mod }[\theta^b]$, for any $0\leq a,b\leq 3$. By reasoning as in Section \ref{trois} we can construct a local diffeomorphism between $\mathfrak{G}_\textsf{f}\times B^4(0,\varepsilon)\simeq \textsf{f}\times B^4(0,\varepsilon)$ and a neighborhood of $\textsf{f}$ in $\mathcal{Y}^{N+1}$. The set $\mathcal{X}^4 =\mathcal{X}$ can thus be endowed with a structure of four-dimensional manifold and $\mathcal{Y}^{N+1}$ with a structure of principal bundle over $\mathcal{X}$ with structure group $\mathfrak{G}_0$.

Again as in Section \ref{trois} a further use of (\ref{ELEYM}b) allows to show the existence of a map $g$ from an open subset of $\mathcal{Y}^{N+1}$ to $\mathfrak{G}_0$ and  of a map $\mathbf{A} = \mathbf{A}_ae^a$, with $\mathbf{A}_a = \mathbf{A}{^i}_a\mathbf{t}_i$,
such that $\theta_z = g^{-1}(z)\mathbf{A}_a(z) e^ag(z) + g^{-1}(z)\hbox{d}g_z$ and, by (\ref{ELEYM}a) $\mathbf{A}_a$ depends only on $x$ in $\mathcal{X}^4$, so that $\mathbf{A} = g\theta g^{-1} - \hbox{d}g\ g^{-1}$ is a 1-form on $\mathcal{X}^4$. We define $\mathbf{F}:= \hbox{d}\mathbf{A} + \frac{1}{2}[\mathbf{A}\wedge \mathbf{A}] = \frac{1}{2}\mathbf{F}_{ab}\theta^a\wedge \theta^b$.

We now let $(S^I_J)_{0\leq I,J\leq N}$ such that
$S^a_b = \delta^a_b$, $S^a_j = S^i_b = 0$ and
$S^i_j\mathbf{t}_i = \hbox{Ad}_g(\mathbf{t}_j)= g\mathbf{t}_jg^{-1}$ (note that $S^K_IS^L_J\textsf{h}_{KL} =\textsf{h}_{IJ}$ since $\textsf{k}$ is invariant by $\hbox{Ad}_g$)  and we introduce
\[
\begin{array}{llll}
 e^I & = & S^I_J\theta^J,\quad & \hbox{ i.e. } e^a:= \theta^a\hbox{ and }
 e^i = S^i_j\theta^j \\
 p_I & = & (S^{-1})^J_I\pi_J & \hbox{ i.e. } p_a:= \pi_a\hbox{ and }
 p_i = (S^{-1})^j_i\pi_j \\
 \omega^{IJ} & = & S^I_KS^J_L\varphi^{KL} - \hbox{d}S^I_KS^J_L\textsf{h}^{KL} \\
 \Omega^{IJ} & = & S^I_KS^J_L\Phi^{KL} 
\end{array}
 \]
We observe that $\mathbf{F}^i = S^i_j\Theta^j$ and $\Omega^{IJ} = \hbox{d}\omega^{IJ} + \omega{^I}_K\wedge \omega^{KJ}$.
We define $\hat{e}^{(N+1)}:= e^0\wedge \cdots \wedge e^N$ and
$\hat{e}^{(N)}_I$, $\hat{e}^{(N-1)}_{IJ}$ and $\hat{e}^{(N-2)}_{IJK}$ in the same way as previously.
We note that, by (\ref{unimodular}), $\hat{e}^{(N+1)} = \hat{\theta}^{(N+1)}$, $\hat{e}^{(N)}_I = (S^{-1})^J_I \hat{\theta}^{(N)}_J$,  and $\hat{e}^{(N-1)}_{IJK} = (S^{-1})^{I'}_I (S^{-1})^{J'}_J (S^{-1})^{K'}_K\hat{\theta}^{(N-1)}_{I'J'K'}$.

By setting $\hbox{d}^\mathbf{A}p_a = \hbox{d}p_a$ and $\hbox{d}^\mathbf{A}p_i = \hbox{d}p_i - c^j_{ki}\mathbf{A}^k\wedge
 p_j$ (\ref{ELEYM}d) translates as
\begin{equation}\label{dpiE}
 \hbox{d}^\mathbf{A} p_I + \frac{1}{2}\hat{e}^{(N-2)}_{IJK}\wedge\Omega^{JK} - \mathbf{F}{^J}_{Ib}\ p{_J}^{jb}\ \hat{e}^{(N)}_j
= \Lambda_0\hat{e}^{(N)}_I
\end{equation}
The computation of $\hbox{d}^\mathbf{A} p_I$ leads to the same result as in the Yang--Mills case (\ref{15abis}), by replacing index $4\leq i\leq N$ by $0\leq I\leq N$ and with the extra simplification that coefficients $p{_I}^{ab}$ vanish.
We observe that $\frac{1}{2}\hat{e}^{(N-2)}_{IJK}\wedge \Omega^{JK} = -\hbox{Ein}(\mathbf{h}){_I}^Je^{(N)}_J$, where 
$\hbox{Ein}(\mathbf{h}){_I}^J = \hbox{Ric}(\mathbf{h}){_I}^J - \hbox{R}(\mathbf{h})\delta{_I}^J$ is the Einstein tensor for the metric $\mathbf{h}$ on $\mathcal{Y}^{N+1}$. Hence, by splitting indices (\ref{dpiE}) thus gives
\begin{equation}\label{decompEinstein}
 \begin{array}{lclc}
  \hbox{Ein}(\mathbf{h}){_a}^b + \Lambda_0\delta{_a}^b 
 & = &  \partial_kp{_a}^{bk} & \hbox{(a)} \\
\hbox{Ein}(\mathbf{h}){_i}^b
 & = & \partial_kp{_i}^{bk} & \hbox{(b)} \\
\hbox{Ein}(\mathbf{h}){_a}^j 
 & = & \partial^\mathbf{A}_cp{_a}^{jc}
+ \partial_kp{_a}^{jk} +
\frac{1}{2} p{_a}^{k\ell}c^j_{k\ell}- \mathbf{F}{^J}_{ac}p{_J}^{jc} & \hbox{(c)}\\
\hbox{Ein}(\mathbf{h}){_i}^j + \Lambda_0\delta{_i}^j
 & = & \partial^\mathbf{A}_cp{_i}^{jc}
+ \partial_kp{_i}^{jk} +
\frac{1}{2}p{_i}^{k\ell}c^j_{k\ell} & \hbox{(d)}
 \end{array}
\end{equation}
where
\begin{equation}\label{17bis}
 \begin{array}{l}
  \partial^{\mathbf{A}}_c p{_a}^{jc} = \partial_cp{_a}^{jc} + \mathbf{A}{^k}_c c^j_{k\ell}p{_a}^{\ell c}  \\
\partial^{\mathbf{A}}_cp{_i}^{jc} = \partial_c p{_i}^{jc} - \mathbf{A}{^k}_c c^\ell_{ki}p{_\ell}^{jc} + \mathbf{A}{^k}_c c^j_{k\ell} p{_i}^{\ell c} 
\end{array}
\end{equation}
The key point now is that, on the one hand, $\hbox{Ein}(\mathbf{h}){_I}^J$ is \emph{constant} on any fiber $\textsf{f} = P^{-1}(x)$ and, on the other hand, since $\textsf{f}$ is compact without boundary,
\[
 \int_\textsf{f}\partial_kp{_I}^{bk}\hat{e}^{(N+1)}
 = \int_\textsf{f}\hbox{d}(p{_I}^{bk}\hat{e}^{(N)}_k) = 0
\]
Hence by the same reasoning as in (\ref{cancel}) we deduce from (\ref{decompEinstein} a,b) that 
\begin{equation}\label{einsteinzero}
  \hbox{Ein}(\mathbf{h}){_a}^b + \Lambda_0\delta{_a}^b =
\hbox{Ein}(\mathbf{h}){_i}^b = 0
\end{equation}
(We also deduce \emph{a posteriori} that $\hbox{Ein}(\mathbf{h}){_a}^j = 0$ because of the symmetry of the Einstein tensor.)

The remaining task is to recognize that (\ref{einsteinzero}) is equivalent to the Einstein--Yang--Mills system on $\mathcal{X}^4$ equipped with the metric $\mathbf{g} = \eta_{ab}e^a\otimes e^b$ and the connection $\mathbf{A}$. This computation is performed (in the tensorial language) in \cite{kerner} and \cite{bleeker}. We present it here using the vierbein formalism.

The spin connection coefficients on $\mathcal{X}^4$ for the metric $\mathbf{g}$ are $\gamma{^a}_b = \gamma{^a}_{bc}e^c$ with
\[
 \gamma{^a}_{bc} = \frac{1}{2}(\Theta{^a}_{bc} - \eta^{ad}\eta_{be}\Theta{^e}_{dc} - \eta^{ad}\eta_{ce}\Theta{^e}_{db})
\]
The relation between the $\gamma{^a}_b$'s and the spin connection $\omega^{IJ} = \omega{^I}_K\textsf{h}^{KJ}$ of $\mathbf{h}$ on $\mathcal{Y}^{N+1}$ is given by
\begin{equation}
 \begin{array}{cclccl}
  \omega{^a}_b & = & \gamma{^a}_b - \frac{1}{2}\textsf{k}_{ij}\eta^{ac}\mathbf{F}{^j}_{cb}e^i \quad &
  \omega{^a}_i & = & \frac{1}{2}\textsf{k}_{ij}\eta^{ac}\mathbf{F}{^j}_{bc}e^b \\
  \omega{^i}_a & = & \frac{1}{2} \mathbf{F}{^i}_{ab}e^b &
  \omega{^i}_j & = & \frac{1}{2}c^i_{jk}(e^k-2\mathbf{A}^k)
 \end{array}
\end{equation}
This leads to the expression of the Einstein tensor of $\mathbf{h}$ in terms of the Einstein tensor of $\mathbf{g}$:
\begin{equation}
 \begin{array}{rcl}
  \hbox{Ein}(\mathbf{h}){_a}^b & = & \hbox{Ein}(\mathbf{g}){_a}^b
  - \frac{1}{2}\mathbf{F}{^i}_{ac}\mathbf{F}{_i}^{bc} + \frac{1}{4}
(\|\mathbf{F}\|^2+\langle \textsf{B}, \textsf{k}\rangle)\delta{_a}^b \\
 \hbox{Ein}(\mathbf{h}){_i}^a & = &
\frac{1}{2}\partial^{\gamma,\mathbf{A}}_b\mathbf{F}{_i}^{ab}
\\ \hbox{Ein}(\mathbf{h}){_i}^j & = &
\frac{1}{4}\mathbf{F}{_i}^{ab} \mathbf{F}{^j}_{ab} - \frac{1}{4}c^k_{i\ell}c^j_{km}\textsf{k}^{\ell m}
+ \frac{1}{4}
(\|\mathbf{F}\|^2+\langle \textsf{B}, \textsf{k}\rangle)\delta{_i}^j
 \end{array}
\end{equation}
where $\|\mathbf{F}\|^2 = \frac{1}{2}\mathbf{F}{_i}^{ab}\mathbf{F}{^i}_{ab}$, $\partial^{\gamma,\mathbf{A}}_b\mathbf{F}{_i}^{ab} = \partial_b\mathbf{F}{_i}^{ab} - c^j_{ki}\mathbf{A}{^k}_b\mathbf{F}{_j}^{ab} + \mathbf{F}{_i}^{cb}\gamma{^a}_{cb} + \gamma{^c}_{bc}\mathbf{F}{_i}^{ab}$ and\footnote{In the version published in Letters in Mathematical Physics, $\langle \textsf{B}, \textsf{k}\rangle$ is denoted by $|c|^2$.} $\langle \textsf{B}, \textsf{k}\rangle := \frac{1}{2}c^i_{\ell k}c^\ell_{ij}\textsf{k}^{jk}$.
Hence, (\ref{einsteinzero}) translates as
\begin{equation}\label{einsteinym}
 \begin{array}{rcl}
  \hbox{Ein}(\mathbf{g}){_a}^b + \Lambda\delta{_a}^b
  & = & \frac{1}{2}(\mathbf{F}{^i}_{ac}\mathbf{F}{_i}^{bc} - \frac{1}{2}
\|\mathbf{F}\|^2\delta{_a}^b) \\
\partial^{\gamma,\mathbf{A}}_b\mathbf{F}{_i}^{ab} & = & 0
 \end{array}
\end{equation}
where $\Lambda := \Lambda_0 + \frac{1}{4}\langle \textsf{B}, \textsf{k}\rangle$.
By choosing the signature of $\eta_{ab}$ to be $(-,+,+,+)$
the first equation reads as the Einstein equation in the presence of a cosmological constant $\Lambda$ and the stress-energy tensor of the gauge fields, and the second one is the Yang--Mills equation on $(\mathcal{X}^4,\mathbf{g})$.

Note that $\langle \textsf{B}, \textsf{k}\rangle$ is equal to $\frac{1}{2}B_{jk}\textsf{k}^{jk}$, where $B_{jk}$ is the Killing form of $\mathfrak{g}$. In particular if $\mathfrak{G}$ is compact semi-simple (which is the case e.g. for $SU(k)$ or $SU(2)\times SU(3)$) $B_{jk}$ is negative definite. Since on the other hand $\textsf{k}_{jk}$ must be positive definite (in order to ensure that the energy of the gauge fields be nonnegative), this implies $\langle \textsf{B}, \textsf{k}\rangle<0$, i.e. $\Lambda < \Lambda_0$.\\

\subsection{Gauge symmetries}
The action (\ref{actionfinale}) and the constraint (\ref{constraint}) are invariant by orientation preserving diffeomorphisms $T:\mathcal{Y}^{N+1}\longrightarrow \mathcal{Y}^{N+1}$ acting on fields through pullback $(\theta,\varphi,\pi)\longmapsto (T^*\theta,T^*\varphi,T^*\pi)$. They are also invariant through the transformation $(\theta^I,\varphi^{IJ},\pi_J)\longmapsto (S^I_J\theta^J,S^I_KS^J_L\varphi^{KL},(S^{-1})^J_I\pi_J)$ where $(S^I_J)$ is the matrix of $\hbox{Ad}_g$, for some $g\in \mathfrak{G}$ \emph{which is constant} (If $g$ is not constant, the curvature 2-form $\Phi^{IJ}$ does not transform in a tensorial way.)

\section{Einstein--Maxwell model}\label{cinq}
What changes if we replace the compact simply connected group  $\mathfrak{G}$ by $U(1)$ ? 
The action in (\ref{actionfinale}) becomes
\begin{equation}
 \mathcal{A}_0[\theta,\varphi,\pi] = \int_{\mathcal{Y}^5} \frac{1}{2}\hat{\theta}^{(3)}_{IJ}\wedge \Phi^{IJ} + \pi_I\wedge \hbox{d}\theta^I 
\end{equation}
where $\pi_I = -\pi{_I}^a\theta^{(3)}_a$.
Critical points of the action $\mathcal{A}_{\Lambda_0}$ in (\ref{22bis}) satisfy System (\ref{ELEYM}), except that (\ref{ELEYM}b) does not exist. This has no incidence for solving the system obtained by imposing that the restriction of $\theta^a$ on $\textsf{f} $ vanishes since the integral leaves are just curves.
The key point is that $U(1)$ is not simply connected. Hence we cannot conclude that the integral leaves are compact and form a fibration in general. 

However, we shall prove that, if we know that \emph{at least one leaf closes}, then \emph{all leaves close} and are diffeomorphic and we hence get the existence of a space-time $\mathcal{X}^4$ and a fibration of $\mathcal{Y}^5$ over it.

Let $\textsf{Y}$ be the vector field such that $\theta^a(\textsf{Y}) = 0$ for any $0\leq a\leq 3$ and $\theta^4(\textsf{Y}) = 1$.
Assume that some integral leaf $\textsf{f}_0$ closes. It means that there exists a map $u:\mathbb{R}\longrightarrow \mathcal{Y}^5$ which is a solution of $\frac{du}{ds}= \textsf{Y}(u)$, the image of which is $\textsf{f}_0$ and which is periodic, i.e., that there exists $q>0$ such that $u(t+q) = u(t)$. W.l.g. we can assume that $q$ is minimal. For any $\xi = (\xi^a)_{0\leq a\leq 3}$ in $\mathbb{R}^4$, let $\textsf{X}$ be the vector field defined by $\theta^a(\textsf{X}) = \xi^a$ for $0\leq a\leq 3$ and $\theta^4(\textsf{X})=0$. 
Then, $[\textsf{X},\textsf{Y}]=0$ since $\hbox{d}\theta^I(\textsf{X},\textsf{Y}) = 0$ because of (\ref{ELEYM}a).
Consider its flow map $e^{\cdot \textsf{X}}$ and, for a fixed $t\in \mathbb{R}$, the map $v:\mathbb{R}\longrightarrow \mathcal{Y}^5$ defined by $v(s) = e^{t\textsf{X}}(u(s))$. Since $\textsf{X}$ and $\textsf{Y}$ commute $\frac{dv}{ds} = \textsf{Y}(v)$. This shows that $e^{t\textsf{X}}(\textsf{f}_0)$, the image of $\textsf{f}_0$ by $e^{t\textsf{X}}$, is an integral leaf $\textsf{f}$ of $\textsf{X}$. Since $\textsf{f}_0$ is compact and $e^{t\textsf{X}}$ is continuous $\textsf{f}$ is compact and hence does not intersect from $\textsf{f}_0$ if $t$ is sufficiently small.  Since this works for all values of $\xi$ and since $\mathcal{Y}^5$ is connected, this endows $\mathcal{Y}^5$ with a topological bundle structure. By proceeding as in Section \ref{deux} we can achieve the normalization (\ref{fluxMaxwell}) and then obtain a solution of the Einstein--Maxwell system by using arguments from Section \ref{trois} and \ref{quatre}.

By replacing $\mathfrak{G}$ by, for example, $U(1)\times SU(2)\times SU(3)$ the situation is similar: we do not obtain a fibration in general, but \emph{we do} if \emph{at least one leaf} closes.

\section*{Conclusion}
Models presented in Section \ref{deux} and \ref{trois} provide gauge theories which, in addition to the usual gauge symmetries, are invariant by the action of diffeomorphisms which respect a topological fibration, enforcing thus the similarity between gauge theories and general relativity. Models in Section \ref{quatre} and \ref{cinq} unify gravity and gauge theories at the classical level without symmetry hypotheses. All these models are based on the introduction of auxiliary `\emph{ghosts}' fields (but not in the Faddeev--Popov sense), which are not observable at the classical level.
The presence of these `ghosts' fields enforces the existence of a principal bundle structure, through the dynamical equations. The key point is that, when projecting the dynamical equations on the base manifolds, these ghosts fields cancel in the equations.
One may ask whether these 'ghosts' fields could create physically observable phenomena in a quantum version of the theory.

An intriguing fact is that, although invariant by diffeomorphisms our Kaluza--Klein theory is not invariant by standard gauge transformation. This may be connected to the difficulty to lift the theory on the bundle of orthonormal frames over space-time (or its spin cover) in the spirit of the  \emph{group manifold} approach as initiated in \cite{neemanregge,toller} and as done in 
\cite{heleinvey15} and further extended in the presence of spinors in \cite{pierard-dirac}.

\bibliographystyle{plain}
\bibliography{bibneo}

\begin{thebibliography}{10}

\bibitem{appelquist}
T.~Appelquist, A.~Chodos, and P.G.O. Freund.
\newblock {\em Modern Kaluza-Klein theories}.
\newblock Addison-Wesley Pub. Co., Boston, 1987.

\bibitem{bleeker}
D.~Bleeker.
\newblock {\em Gauge Theory and Variational Principles}.
\newblock Addison-Wesley, Mass., 1981.

\bibitem{bcggg}
R.~L. Bryant, S.~S. Chern, R.~B. Gardner, H.~L. Goldschmidt, and P.~A.
  Griffiths.
\newblock {\em Exterior Differential Systems}.
\newblock Springer Verlag, New-York,Paris, 1991.

\bibitem{helein14}
F.~H{\'{e}}lein.
\newblock Multisymplectic formulation of {Y}ang-{M}ills equations and
  {E}hresmann connections.
\newblock {\em Advances in Theoretical and Mathematical Physics},
  19(4):805--835, 2015.
\newblock \href{https://arxiv.org/abs/1406.3641}{arXiv:1406.3641}.

\bibitem{helein2020}
F.~H{\'{e}}lein.
\newblock A variational principle for kaluza{\textendash}klein types theories.
\newblock {\em Advances in Theoretical and Mathematical Physics},
  24(2):305--326, 2020.
\newblock \href{https://arxiv.org/abs/1809.03375}{arXiv:1809.03375}.

\bibitem{heleinvey15}
F.~H{\'{e}}lein and D.~Vey.
\newblock Curved space-times by crystallization of liquid fiber bundles.
\newblock {\em Foundations of Physics}, 47(1):1--41, sep 2016.
\newblock \href{https://arxiv.org/abs/1508.07765}{arXiv:1508.07765}.

\bibitem{jordan}
P.~Jordan.
\newblock Erweiterung der projektiven {R}elativit\"{a}tstheorie.
\newblock {\em Annalen der Physik}, 436(4-5):219--228, 1947.

\bibitem{kaluza}
Th. Kaluza.
\newblock On the unification problem in physics.
\newblock {\em Sitzungsberichte Pruss. Acad. Sci.}, page 966, 1921.
\newblock Reprinted in English in \cite{appelquist} and in
  \href{https://arxiv.org/abs/1803.08616}{arXiv:1803.08616}.

\bibitem{kerner}
R.~Kerner.
\newblock Generalization of the {K}aluza-{K}lein theory for an arbitrary
  non-{A}belian gauge group.
\newblock {\em Ann. Inst. H. Poincar{\'e}}, 9(2):143, 1968.
\newblock in
  \href{http://www.numdam.org/item?id=AIHPA_1968__9_2_143_0}{www.numdam.org}.

\bibitem{klein}
O.~Klein.
\newblock Quantentheorie und f{\"u}nfdimensionale {R}elativit{\"a}tstheorie.
\newblock {\em Zeitschrift f{\"u}r Physik}, 37(12):895--906, dec 1926.

\bibitem{neemanregge}
Y.~Ne'eman and T.~Regge.
\newblock Gauge theory of gravity and supergravity on a group manifold.
\newblock {\em La Rivista Del Nuovo Cimento}, 1(5):1--43, 1978.

\bibitem{pierard-dirac}
Jérémie Pierard~de Maujouy.
\newblock Dirac spinors on generalised frame bundles: a frame bundle
  formulation for einstein-cartan-dirac theory.

\bibitem{thiry}
Y.~Thiry.
\newblock Les {\'e}quations de la th{\'e}orie unitaire de {K}aluza.
\newblock {\em Comptes Rendus Acad. Sci. Paris}, 226(216), 1948.

\bibitem{toller}
M.~Toller.
\newblock Classical field theory in the space of reference frames.
\newblock {\em Il Nuovo Cimento B Series 11}, 44(1):67--98, mar 1978.

\end{thebibliography}

\end{document}